\documentclass[thmsa,11pt]{article}
\usepackage{latexsym}
\usepackage{amsmath,amssymb,amscd}
\usepackage{graphics}
\usepackage{epsfig}
\usepackage{graphicx}
\usepackage{wrapfig}
\usepackage[rflt]{floatflt}
\usepackage{pgf}
\usepackage{tikz}
\usepackage[latin1]{inputenc}

\usetikzlibrary{arrows,automata}
\newtheorem{theorem}{Theorem}[section]

\newtheorem{lemma}{Lemma}[section]
\newtheorem{proposition}{Proposition}[section]
\newtheorem{definition}{Definition}[section]

\newenvironment{proof}{\trivlist\item[\hskip
       \labelsep{\it Proof:\/}]\ignorespaces}{\hfill$\Box$\endtrivlist}

\newcommand{\Go}[1]{{G\"{o}del} }

\begin{document}

\date{}
\title{Bi-modal Gödel logic over [0,1]-valued \\
Kripke frames}
\author{Xavier Caicedo \thanks{{\small Departamento de Matemáticas,
Universidad de los Andes, Bogotá, Colombia;} {\small \texttt{%
xcaicedo@uniandes.edu.co} }} \and Ricardo Oscar Rodríguez$\ $\thanks{{\small %
Departamento de Computación, Fac. Ciencias Exactas y Naturales,} {\small %
Universidad de Buenos Aires, 1428 Buenos Aires, Argentina;\ } {\small
\texttt{ricardo@dc.uba.ar} }}}
\maketitle

\begin{abstract}
We consider the Gödel bi-modal logic determined by fuzzy Kripke models where
both the propositions and the accessibility relation are infinitely valued
over the standard Gödel algebra [0,1] and prove strong completeness of
Fischer Servi intuitionistic modal logic IK plus the prelinearity axiom with
respect to this semantics. We axiomatize also the bi-modal analogues of $T,$
$S4,$ and $S5$ obtained by restricting to models over frames satisfying the
[0,1]-valued versions of the structural properties which characterize these
logics. As application of the completeness theorems we obtain a
representation theorem for bi-modal Gödel algebras.
\end{abstract}

In a previous paper \cite{CaicedoRodriguez}, we have considered a semantics
for Gödel modal logic based on fuzzy Kripke models where both the
propositions and the accessibility relation take values in the standard Gö%
del algebra [0,1], we call these Gödel-Kripke models, and we have\ provided
strongly complete axiomatizations for the uni-modal fragments of this logic
with respect to validity and semantic entailment from countable theories.
The systems $\mathcal{G}_{\square }$ and $\mathcal{G}_{\Diamond }$
axiomatizing the $\square $-fragment and the $\Diamond $-fragment,
respectively, are obtained by adding to Gödel-Dummet\ propositional
calculus\ the following axiom schemes and inference rules:%
\begin{equation*}
\begin{array}{ll}
\mathcal{G}_{\square }\text{:} & \square (\varphi \rightarrow \psi
)\rightarrow (\square \varphi \rightarrow \square \psi ) \\
& \lnot \lnot \square \varphi \rightarrow \square \lnot \lnot \varphi \\
& \emph{From\ }\varphi ,\emph{\ infer\ }\square \varphi \text{ } \\
&
\end{array}%
\text{ \ }%
\begin{array}{ll}
\mathcal{G}_{\Diamond }\text{:} & \Diamond (\varphi \vee \psi )\rightarrow
(\Diamond \varphi \vee \Diamond \psi ) \\
& \Diamond \lnot \lnot \varphi \rightarrow \lnot \lnot \Diamond \varphi \\
& \lnot \Diamond \bot \\
& \emph{From\ }\varphi \rightarrow \psi ,\emph{\ infer\ }\Diamond \varphi
\rightarrow \Diamond \psi%
\end{array}%
\end{equation*}%
These logics diverge substantially in their model theoretic properties.
Thus,\ $\mathcal{G}_{\square }$ does not have the finite model property
while $\mathcal{G}_{\Diamond }$ does, and the first logic is characterized
by models with \{0,1\}-valued accessibility relation (accessibility-crisp
models) while the second does not. Similar results were obtained for the
uni-modal Gödel analogues of the classical modal logics $T$ and $S4$
determined by Gödel-Kripke models over frames satisfying, respectively, the
[0,1]-valued version of reflexivity, or reflexivity and transitivity. The
axiomatization of the uni-modal Gödel analogues of $S5$ remains open.

It is the main purpose of this paper to show that the full bi-modal logic
based in Gödel-Kripke models is axiomatized by the system $\mathcal{G}%
_{\square \Diamond }$ which results of adding to the union of $\mathcal{G}%
_{\square }$ and $\mathcal{G}_{\Diamond }$ Fischer-Servi´s connecting axioms
\cite{FisherServi84}:\medskip

$%
\begin{array}{l}
\Diamond (\varphi \rightarrow \psi )\rightarrow (\Box \varphi \rightarrow
\Diamond \psi ) \\
(\Diamond \varphi \rightarrow \Box \psi )\rightarrow \Box (\varphi
\rightarrow \psi ),%
\end{array}%
$\medskip

\noindent and to extend this completeness result to the bi-modal Gödel
analogues of classical $T,$ $S4,$ and $S5.$

The many valued Kripke interpretation of bi-modal logic utilized in this
paper was proposed originally by Fitting \cite{Fitting1}, \cite{Fitting2}
with a complete Heyting algebra as algebra of truth values, and he gave a
complete axiomatization assuming the algebra was finite and the language had
constants for all the truth values. In \cite{Koutras} and \cite{Koutras2}
transformations and characterization of frame properties for these models
are given in the general case. Bou, Esteva, and Godo \cite{BouEstevaGodo}
have proposed to utilize this kind of interpretation for general algebras in
the study of fuzzy modal logics.

Our method of proof do not seem to extend easily, however, to algebras
distinct from the Gödel algebra [0,1] and we do not know any other
completeness result for this type of semantics for a fixed algebra $H$,
except Fitting´s quoted above and Metcalfe \& Olivetti completeness of a
natural deduction system for the $\square $-fragment of our Gödel-modal
logic \cite{Metcalfe}.

$\mathcal{G}_{\square \Diamond }$ may be shown deductively equivalent to the
system $IK$ introduced by Fischer-Servi \cite{FisherServi84} as the natural
intuitionistic counterpart of classical modal logic, plus the prelinearity
axiom: $(\varphi \rightarrow \psi )\vee (\psi \rightarrow \varphi ).$
Similarly, the Gödel analogue of bi-modal $S5$ results equivalent to the
system $MIPC$ of Prior \cite{Prior} plus prelinearity.

$IK$ and its extensions have been extensively studied, either by means of
classical Kripke models for intuitionism equipped with extra relations
commuting with the order to interpret the modal operators. (\cite{Simpson},
\cite{Ono77}, \cite{OnoSuzuki}, \cite{Dosen}, \cite{Wolter}, \cite%
{WolterZakha}, \cite{Grefe98}, \cite{Celani}, \cite{Davoren}), or by means
of algebraic interpretations, specially in the case of $MIPC,$ known to be
complete for values in monadic Heyting algebras (\cite{Bull}, \cite{Ono77},
\cite{FisherServi78}, \cite{Guram1}, \cite{Guram2}) A major result is that
under these semantics $IK$ and $MIPC$ enjoy the finite model property.

Clearly, $\mathcal{G}_{\square \Diamond }$ and its modal extensions inherit
these semantics\ by asking the multirelational Kripke frames to be linearly
ordered or the algebras to be Gödel algebras, but these alternative
interpretations do not have the standard character of Gödel-Kripke semantics
relevant to fuzzy logic, nor seem our results reducible to their properties.
For example, the formula $\square \lnot \lnot \theta \rightarrow \lnot \lnot
\square \theta $ has finite counter-models in them but not in Gödel-Kripke
semantics. We discuss briefly at the end of the paper an embedding of our
semantics into algebraic semantics and utilize our completeness theorem to
show a representation theorem for countable bi-modal Gödel algebras.

\section{Gödel Kripke models}

The language\textit{\ }$\mathcal{L}_{\square \Diamond }(Var)$ of
propositional \emph{bi-modal logic} is built from a set $Var$ of
propositional variables, connectives symbols $\vee ,\wedge ,\rightarrow
,\bot ,$ and the modal operators symbols $\square $ and $\Diamond $. Other
connectives are defined as usual: $\top :=\varphi \rightarrow \varphi ,$ $%
\lnot \varphi :=\varphi \rightarrow \bot ,$ $\varphi \longleftrightarrow
\psi :=(\varphi \rightarrow \psi )\wedge (\psi \rightarrow \varphi ).$ We
will write $\mathcal{L}_{\square \Diamond }$ if the set $Var$ is understood.

Recall that a \emph{linear Heyting algebra}, or \emph{Gödel algebra} in the
fuzzy literature, is a Heyting algebra satisfying the identity $%
(x\Rightarrow y)\curlyvee (y\Rightarrow x)=1.$ The variety of these algebras
is generated by the \emph{standard Gödel algebra }$[0,1]$, the ordered
interval with its unique Heyting algebra structure. Let the symbols $\cdot $%
, $\Rightarrow ,\curlyvee ,$ and denote, respectively, the meet, residuum
(implication), and join operations of $[0,1]$.\footnote{%
The join operation is definable in Gödel algebras as\ $x\curlyvee
y=((x\Rightarrow y)\Rightarrow y))\cdot ((y\Rightarrow x)\Rightarrow x))$}

\begin{definition}
\label{models} \emph{A} Gödel-Kripke model ($GK$-model)\ \emph{will be a
structure }$M=\langle W,S,e\rangle $\emph{\ where }$W$\emph{\ is a non-empty
set of objects that we call }worlds of\emph{\ }$M,$\emph{\ and }$S:W\times
W\rightarrow \lbrack 0,1],$ $e:W\times Var\rightarrow \lbrack 0,1]$\emph{\
are arbitrary\ functions. The pair }$\langle W,S\rangle $\emph{\ will be
called a }$GK$-frame.
\end{definition}

\noindent The function $e:W\times Var\rightarrow \lbrack 0,1]$\ associates
to each world $x$ a valuation $e(x,-):Var\rightarrow \lbrack 0,1]$ which
extends to $e(x,-):\mathcal{L}_{\square \Diamond }(Var)\rightarrow \lbrack
0,1]$\ by defining inductively on the construction of the formulas (we
utilize the same symbol $e$ to name the extension):

\medskip

$e(x,\bot ):=0$

$e(x,\varphi \wedge \psi ):=e(x,\varphi )\cdot e(x,\psi )$

$e(x,\varphi \vee \psi ):=e(x,\varphi )\curlyvee e(x,\psi )$

$e(x,\varphi \rightarrow \psi ):=e(x,\varphi )\Rightarrow e(x,\psi )$

$e(x,\bot ):=0$

$e(x,\Box \varphi ):=\inf_{y\in W}\{Sxy\Rightarrow e(y,\varphi )\}$

$e(x,\Diamond \varphi ):=\sup_{y\in W}\{Sxy\cdot e(y,\varphi )\}$.

\medskip

\noindent Truth, validity and entailment are defined as follows for $\varphi
\in \mathcal{L}_{\square \Diamond },$ $T\subseteq \mathcal{L}_{\square
\Diamond }$:\medskip

- $\varphi $ is \emph{true in} $M$ \emph{at} $x,$ written $M\models
_{x}\varphi $, if $e(x,\varphi )=1.$

- $\varphi $ is \emph{valid in} $M,$ written $M\models \varphi $, if $%
M\models _{x}\varphi $ at any world $x$ of $M.$

- $\varphi $ is $GK$-\emph{valid}$,$ written $\models _{GK}\varphi $, if $%
M\models \varphi $ for any $GK$-model $M$.

- $T\models _{GK}\varphi $ if and only if for any $GK$-model $M$ and any
world $x$ in $M$ :

$\ \ M\models _{x}\theta $ for all $\theta \in T$ implies $M\models
_{x}\varphi .$\medskip

It is routine to verify that all axiom schemes corresponding to identities
satisfied by $[0,1];$ that is, the laws of Gödel-Dummet logic, are $GK$%
-valid. In addition

\begin{proposition}
\label{soundness}The following schemes are\textbf{\ }$GK$-valid:
\end{proposition}

\noindent $%
\begin{array}{lll}
\text{K}_{\square } &  & \square (\varphi \rightarrow \psi )\rightarrow
(\square \varphi \rightarrow \square \psi ) \\
\text{K}_{\Diamond } &  & \Diamond (\varphi \vee \psi )\rightarrow (\Diamond
\varphi \vee \Diamond \psi ) \\
\text{F}_{\Diamond } &  & \lnot \Diamond \bot \\
\text{FS1} &  & \Diamond (\varphi \rightarrow \psi )\rightarrow (\Box
\varphi \rightarrow \Diamond \psi ) \\
\text{FS2} &  & (\Diamond \varphi \rightarrow \Box \psi )\rightarrow \Box
(\varphi \rightarrow \psi ).%
\end{array}%
$

\proof Let $M=\langle W,S,e\rangle $ be a $GK$-model. (K$_{\square })$ By
definition and properties of the residuum, $e(x,\square (\varphi \rightarrow
\psi ))\cdot e(x,\square \varphi )$ $\leq (Sxy\Rightarrow (e(y,\varphi
)\Rightarrow e(y,\psi ))\cdot (Sxy\Rightarrow e(y,\varphi ))$ $\leq
(Sxy\Rightarrow e(y,\psi ))$ for any $y\in W$. Taking the meet over $y$ in
the last expression: $e(x,\square (\varphi \rightarrow \psi ))\cdot
e(x,\square \varphi )$ $\leq e(x,\square \psi ),$ hence\ $e(x,\square
(\varphi \rightarrow \psi ))$ $\leq e(x,\square \varphi \rightarrow \square
\psi )$. (K$_{\Diamond })$\ By distributivity and properties of the join: $%
e(\Diamond (\varphi \vee \psi ))=\sup_{y}\{Sxy\cdot (e(y,\varphi )\curlyvee
e(y,\psi ))\}=$ $\sup_{y}\{Sxy\cdot e(y,\varphi )\}\curlyvee
\sup_{y}\{Sxy\cdot e(y,\psi )\}.$ ($\mathbf{F}_{\Diamond })$ $e(x,\Diamond
\bot )=\sup_{y}\{Sxy\cdot 0\}=0.$ (FS1) $Sxy\cdot e(x,\Box \varphi )\cdot
e(y,\varphi \rightarrow \psi )\leq Sxy\cdot (Sxy\Rightarrow e(y,\varphi
))\cdot (e(y,\varphi )\Rightarrow e(y,\psi ))$ $\leq Sxy\cdot e(y,\psi
))\leq e(x,\Diamond \psi ).$ Therefore, $Sxy\cdot e(y,\varphi \rightarrow
\psi )\leq (e(x,\Box \varphi )\Rightarrow e(x,\Diamond \psi )),$ and taking
the join over $y$ in the left hand side, we have $e(x,\Diamond (\varphi
\rightarrow \psi ))\leq e(x,\Box \varphi \rightarrow \Diamond \psi ).$ (FS2)
$e(x,\Diamond \varphi \rightarrow \Box \psi )\leq $ $[Sxy\cdot e(y,\varphi
)\Rightarrow (Sxy\Rightarrow e(y,\psi ))]$ $=[Sxy\cdot e(y,\varphi
)\Rightarrow e(y,\psi )]$ $=(Sxy\Rightarrow (e,y\rightarrow \psi )).$ $%
\blacksquare $

\emph{\medskip }

\textbf{Remark}. Changing the algebra $[0,1]$ to a complete Heyting algebra $%
H$ in the above definitions we have Kripke models valued in a $H$ ($HK$%
-models) and the corresponding notion of $HK$-validity. Then all laws of the
intermediate logic determined by $H$ are $HK$-valid, and also Proposition %
\ref{soundness} holds.

\section{A deductive calculus}

Let $\mathcal{G}$ be some axiomatic version of Gödel-Dummet propositional
calculus;\ that is, Heyting calculus plus the prelinearity axiom $(\varphi
\rightarrow \psi )\rightarrow (\psi \rightarrow \varphi )$, and let $\vdash
_{\mathcal{G}}$ denote deduction in this logic. Let $\mathcal{L}(X)$ denote
the set of formulas built by means of the connectives $\wedge ,\rightarrow ,$
and $\bot ,$ from a given set $X.$ For simplicity, the extension of a
valuation $v:X\rightarrow \lbrack 0,1]$ to $\mathcal{L}(X) $ according to
the Heyting interpretation of the connectives will be denoted also $v.$It is
well known that this system is complete for validity with respect to these
valuations and the distinguished value 1. We will utilize the fact that it
is actually sound and complete in the following strong way (see \cite%
{CaicedoRodriguez}):

\begin{proposition}
{\label{ordersoundness} i)\ }If $T\cup \{\varphi \}\subseteq \mathcal{L}(X),$
then $T\vdash _{\mathcal{G}}\varphi $ implies $\inf v(T)\leq v(\varphi )$
for any valuation $v:X\rightarrow \lbrack 0,1]$. ii)\ If $T$ is countable,
and $T\nvdash _{\mathcal{G}}\varphi _{i_{1}}\vee ..\vee \varphi _{i_{1}}$
for each finite subset of a countable family $\{\varphi _{i}\}_{i}$ there is
a valuation $v:L\rightarrow \lbrack 0,1]$\ such that $v(\theta )=1$\ for all
$\alpha \in T$\ and $v(\varphi _{i})<1$ for all $i.$
\end{proposition}

For an example that completeness for [0,1]-valued entailment can not be
extended to uncountable theories see Section 3 in \cite{CaicedoRodriguez}
and also Proposition \ref{uncountable} in this paper.

\begin{definition}
$\mathcal{G}_{\square \Diamond }$ \emph{is the deductive calculus obtained
by adding to }$\mathcal{G}$\emph{\ the schemes }K$_{\square },$ K$_{\Diamond
},$ F$_{\Diamond },$ FS1, FS2 of Proposition \ref{soundness} and the
inference rules:
\end{definition}

\noindent $%
\begin{array}{lll}
\text{NR}_{\square } &  & \emph{From\ }\varphi \emph{\ infer}\text{ }\square
\varphi \text{ } \\
\text{RN}_{\Diamond } &  & \emph{From\ }\varphi \rightarrow \psi \emph{\
infer\ }\Diamond \varphi \rightarrow \Diamond \psi \text{.}%
\end{array}%
$\medskip

\noindent Proofs with assumptions are allowed with the restriction that NR$%
_{\square }$ and RN$_{\Diamond }$ may be applied only when the premise is a
theorem. Let $\vdash _{\mathcal{G}_{\square \Diamond }}$denote deduction in
this system.\medskip

The restriction on the application of the rules allows the following
convenient reduction (see \cite{CaicedoRodriguez}).\medskip

\begin{lemma}
{\label{reduction} }Let $Th\mathcal{G}_{\square \Diamond }$ be the set of
theorems of $\mathcal{G}_{\square \Diamond }$ with no assumptions, then for
any theory $T$ and formula $\varphi $ in $\mathcal{L}_{\Box \Diamond
}:T\vdash _{\mathcal{G}_{\Box \Diamond }}\varphi $\ if\ and\ only\ if\ $%
T\cup Th\mathcal{G}_{\square \Diamond }\vdash _{\mathcal{G}}\varphi $.
\end{lemma}

\noindent and the Deduction Theorem:.

\begin{lemma}
$T,\psi \vdash _{\mathcal{G}_{\Box \Diamond }}\varphi $ implies \ $T\vdash _{%
\mathcal{G}_{\Box \Diamond }}\psi \rightarrow \varphi .$
\end{lemma}

The following are theorems of $\mathcal{G}_{\Box \Diamond }$. The first one
is given as an axiom in Fitting \cite{Fitting1}, the next two show that $%
\mathcal{G}_{\square \Diamond }$ is just the union of $\mathcal{G}_{\square
} $, $\mathcal{G}_{\Diamond }$ plus the Fischer Servi axioms, and the last
one will be useful in our completeness proof. \medskip

\noindent $%
\begin{array}{ll}
\text{T1.} & \lnot \Diamond \theta \longleftrightarrow \square \lnot \theta
\\
\text{T2}. & \lnot \lnot \square \theta \rightarrow \square \lnot \lnot
\theta \\
\text{T3}. & \Diamond \lnot \lnot \varphi \rightarrow \lnot \lnot \Diamond
\varphi \\
\text{T4.} & (\square \varphi \rightarrow \Diamond \psi )\vee \square
((\varphi \rightarrow \psi )\rightarrow \psi )%
\end{array}%
$\medskip

\noindent To see this, we write temporarily $\vdash $ for $\vdash _{\mathcal{%
G}_{\Box \Diamond }}$, then

(T1) \ $\lnot \Diamond \theta \vdash (\Diamond \theta \rightarrow \square
\bot )\vdash \square (\theta \rightarrow \bot )\ $by Heyting calculus and
FS2. Similarly, $\Diamond \theta \vdash \Diamond (\lnot \theta \rightarrow
\bot )\vdash \square \lnot \theta \rightarrow \Diamond \bot \vdash \lnot
\square \lnot \theta \ $by Heyting calculus, RN$_{\Diamond },$ and FS2;
hence, $\square \lnot \theta \vdash \lnot \Diamond \theta $.

(T2) $(\square \varphi \rightarrow \bot )\rightarrow \bot \vdash (\square
\varphi \rightarrow \Diamond \bot )\rightarrow \bot \vdash \Diamond (\varphi
\rightarrow \bot )\rightarrow \square \bot \vdash \square ((\varphi
\rightarrow \bot )\rightarrow \bot )$ by F$_{\Diamond }$, FS2, and FS1.

(T3) \ From FS1,\ $\vdash \Diamond (\lnot \varphi \rightarrow \bot
)\rightarrow (\Box \lnot \varphi \rightarrow \Diamond \bot );$ that is, $%
\vdash \Diamond (\lnot \lnot \varphi )\rightarrow (\lnot \Diamond \varphi
\rightarrow \bot )$ by T1 and F$_{\Diamond }.$

(T4) By prelinearity:\ $\vdash (\square \varphi \rightarrow \Diamond
(\varphi \rightarrow \psi ))\vee (\Diamond (\varphi \rightarrow \psi
)\rightarrow \square \varphi ),$ but $\square \varphi \rightarrow \Diamond
(\varphi \rightarrow \psi )$ $\vdash \square \varphi \rightarrow (\square
\varphi \rightarrow \Diamond \psi )$ $\vdash \square \varphi \rightarrow
\Diamond \psi $ by FS1; moreover, $\Diamond (\varphi \rightarrow \psi
)\rightarrow \square \varphi $ $\vdash \square ((\varphi \rightarrow \psi
)\rightarrow \varphi )$ $\vdash \square ((\varphi \rightarrow \psi
)\rightarrow \psi )$ by FS2, Heyting calculus and RN$_{\square }.$

\begin{theorem}
\label{strong soundness}(Soundness) $T\vdash _{\mathcal{G}_{\Box \Diamond
}}\varphi $ implies $T\models _{GK}\varphi $.
\end{theorem}

\proof Clearly, the Modus Ponens rule preserves truth at every world of any $%
GK$-model $M$. Moreover, $M\models \varphi $ implies $M\models \square
\varphi ,$ trivially, and $M\models \varphi \rightarrow \psi $ implies $%
M\models \Diamond \varphi \rightarrow \Diamond \psi $ because if $%
e(x,\varphi \rightarrow \psi )=1$ for all $x$ then $Sxy\cdot e(y,\varphi
)\leq Sxy\cdot e(y,\psi )\leq e(x,\Diamond \psi )$ for all $x,y,$ and taking
the join in the left, $e(x,\Diamond \varphi )\leq e(x,\Diamond \psi ).$ The
rest follows from Proposition \ref{soundness}. $\blacksquare $\medskip

It is easy to provide counterexamples to the validity of $\lnot \square
\lnot \theta \rightarrow \Diamond \theta $ and $\lnot \Diamond \lnot \theta
\rightarrow \square \theta ,$ thus the modal operators are not
interdefinable in $\mathcal{G}_{\square \Diamond }$ in the classical manner.
In fact, they are not interdefinable in any way. For example,\ the formula $%
\square \lnot \lnot \theta \rightarrow \lnot \lnot \square \theta $ is not
expressible in terms of $\Diamond $ alone because the $\Diamond $-fragment
has the finite model property with respect to the number of worlds while
this formula has not finite counterexamples as shown in \cite%
{CaicedoRodriguez}.\medskip

\noindent \textbf{Remark}. $\mathcal{G}_{\square \Diamond }$ may be seen
deductively equivalent to Fischer-Servi system $IK$ (cf. \cite{FisherServi84}%
) plus the prelinearity axiom, replacing K$_{\square }$ with the axiom $%
\square (\varphi \wedge \psi )\rightarrow (\square \varphi \wedge \square
\psi )$ and the rule and NR$_{\square }$ with $\emph{From\ }\varphi
\rightarrow \psi \emph{\emph{\ }\ infer}$ $\square \varphi \rightarrow $ $%
\square \psi .$ Actually, T1, T2, T3 are theorems of $IK,$ and we have that $%
T\vdash _{IK}\varphi $ implies $T\models _{HK}\varphi $ for any complete
Heyting algebra $H.$

\section{Completeness}

To prove \emph{strong completeness} of $\mathcal{G}_{\square \Diamond }$
with respect to entailment from countable theories in Gödel-Kripke
semantics, our strategy is to show this for finite theories first, and then
utilize a first order compactness argument to lift it to countable theories.
To show weak completeness we define for each finite \emph{fragment} $%
F\subseteq \mathcal{L}_{\square \Diamond }$ (that is, a subset closed under
subformulas and containing the formula $\bot )$ a canonical model.

Denote by $\square \mathcal{L}_{\square \Diamond }\ $and $\Diamond \mathcal{L%
}_{\square \Diamond }$ the sets of formulas in $\mathcal{L}_{\square
\Diamond }$ starting with $\square $ and $\Diamond ,$ respectively, and set $%
X:=\square \mathcal{L}_{\square \Diamond }\cup \Diamond \mathcal{L}_{\square
\Diamond },$ then clearly $\mathcal{L}_{\square \Diamond }(Var)=\mathcal{L(}%
Var\cup X)$. Recall that $Th\mathcal{G}_{\square \Diamond }$ denotes the set
of theorems of $\mathcal{G}_{\square \Diamond }.$

The \emph{canonical model }$M_{F}=(W,S^{F},e^{F})$ is defined as follows.

\begin{description}
\item[$W$:] is the set of valuations $v:Var\cup X\rightarrow \lbrack 0,1]$
such that $v(Th\mathcal{G}_{\square \Diamond })=1$ when $Th\mathcal{G}%
_{\square \Diamond }$ is considered as a subset of $\mathcal{L(}Var\cup X)$.

\item[$S^{F}$:] $S^{F}vw=\inf_{\psi \in F}\{(v(\Box \psi )\rightarrow w(\psi
))\cdot (w(\psi )\rightarrow v(\Diamond \psi ))\}.$

\item[$e^{F}$:] $e^{F}(v,p)=v(p)$ for any $p\in Var$.
\end{description}

Weak completeness will follow from the following lemma which unfortunately
has a rather involved proof.

\begin{lemma}
{\label{equation-joint}}. $e^{F}(v,\varphi )=v(\varphi )$ for any $\varphi
\in F$ and any $v\in W$. $.$
\end{lemma}

\proof We prove this by induction in the complexity of the formulas in $F,$
now considered a subset of $\mathcal{L}_{\square \Diamond }(Var)$. For $\bot
$ and the propositional variables in $F$ the equation holds by definition.
The only non trivial inductive steps are:\ $e^{F}(v,\Box \varphi )=v(\Box
\varphi )$ and $e^{F}(v,\Diamond \varphi )=v(\Diamond \varphi )$ for $\Box
\varphi ,\Diamond \varphi \in F.$ By the inductive hypothesis we may assume
that $e^{F}(v^{\prime },\varphi )=v^{\prime }(\varphi )$ for every $%
v^{\prime }\in W,$ thus we must prove
\begin{equation}
\inf_{v^{\prime }\in W}\{S^{F}vv^{\prime }\Rightarrow v^{\prime }(\varphi
)\}=v(\Box \varphi ))  \label{box}
\end{equation}%
\begin{equation}
\sup_{v^{\prime }\in W}\{S^{F}vv^{\prime }\cdot v^{\prime }(\varphi
)\}=v(\Diamond \varphi ))  \label{Diam}
\end{equation}%
By definition, $S^{F}vv^{\prime }\leq (v(\Box \varphi )\Rightarrow v^{\prime
}(\varphi ))$ and $S^{F}vv^{\prime }\leq (v^{\prime }(\varphi )\Rightarrow
v(\Diamond \varphi ))$ for any $\varphi \in F$ and $v^{\prime }\in W;$
therefore, $v(\Box \varphi )\leq (S^{F}vv^{\prime }\Rightarrow v^{\prime
}(\varphi ))$ and $S^{F}vv^{\prime }\cdot v^{\prime }(\varphi )\leq
v(\Diamond \varphi ).$ Taking the meet over $v^{\prime }$ in the first
inequality and the join in the second,
\begin{equation*}
v(\Box \varphi )\leq \inf_{v^{\prime }\in W}\{S^{F}vv^{\prime }\Rightarrow
v^{\prime }(\varphi )\}\text{, \ }\sup_{v^{\prime }\in W}\{S^{F}vv^{\prime
}\cdot v^{\prime }(\varphi )\}\leq v(\Diamond \varphi ).
\end{equation*}%
Hence, if $v(\Box \varphi )=1$ and $v(\Diamond \varphi )=0$ we obtain (\ref%
{box}) and (\ref{Diam}), respectively. Therefore, it remains only to prove
the next two claims for $\Box \varphi ,\Diamond \varphi \in F$.\medskip

\noindent \textbf{Claim 1}. \emph{If }$v(\Box \varphi )=\alpha <1$\emph{\
and }$\varepsilon >0$ \emph{there exists a valuation }$w\in W$\emph{\ such
that }$S^{F}vw>w(\varphi )$\emph{\ and }$w(\varphi )<\alpha +\varepsilon $.
\emph{That is,} $(S^{F}vw\Rightarrow w(\varphi ))<\alpha +\varepsilon $%
.\medskip

\noindent \textbf{Claim 2. }\emph{If }$v(\Diamond \varphi )=\alpha >0$ \emph{%
then for any }$\varepsilon >0$ \emph{there exists }$w\in W$\emph{\ such that
}$S^{F}vw\cdot w(\varphi )\geq \alpha -\varepsilon .$\medskip

\noindent \emph{Proof of Claim1}. Assume $v(\Box \varphi )=\alpha <1$ and
define (all formulas involved belonging to\ $\mathcal{L}_{\square \Diamond
}(Var)$)
\begin{equation*}
\begin{array}{ll}
\Gamma _{\varphi ,v}= & \{\theta :v(\Box \theta )>\alpha \}\cup \{\theta
_{1}\rightarrow \theta _{2}:v(\Diamond \theta _{1})\leq v(\square \theta
_{2})\} \\
& \cup \{(\theta _{2}\rightarrow \theta _{1})\rightarrow \theta
_{1}:v(\Diamond \theta _{1})<v(\square \theta _{2})\}.%
\end{array}%
\end{equation*}%
Then we have $v(\square \xi )>\alpha $ for each $\xi \in \Gamma _{\varphi
,v}.$ For the first set of formulas by construction. For the second because $%
v(\square (\theta _{1}\rightarrow \theta _{2}))\geq v(\Diamond \theta
_{1}\rightarrow \square \theta _{2})=1$ by FS$2$. For the third, because $%
v(\square \theta _{2}\rightarrow \Diamond \theta _{1})<1$ and thus $%
v(\square ((\theta _{2}\rightarrow \theta _{1})\rightarrow \theta _{1}))=1$
by $T4$. This implies that
\begin{equation*}
\Gamma _{\varphi ,v}\not\vdash _{\mathcal{G}_{\square \Diamond }}\varphi .
\end{equation*}%
Otherwise it would exist $\xi _{1},\ldots ,\xi _{k}\in \Gamma _{\varphi ,v}$
such that: $\xi _{1},\ldots ,\xi _{k}\vdash _{\mathcal{G}_{\square \Diamond
}}\varphi .$ Hence, $\Box \xi _{1},\ldots ,\Box \xi _{k}\vdash _{\mathcal{G}%
_{\square \Diamond }}\Box \varphi $, that is $\Box \xi _{1},\ldots ,\Box \xi
_{k},Th\mathcal{G}_{\square \Diamond }\vdash \Box \varphi $ by Lemma \ref%
{reduction} and thus by Proposition \ref{ordersoundness} (i)%
\begin{equation*}
\alpha <\inf v(\{\Box \xi _{1},\ldots ,v(\Box \xi _{k})\}\cup Th\mathcal{G}%
_{\square \Diamond })\leq v(\square \varphi )=\alpha ,
\end{equation*}%
a contradiction. Therefore, there exists by Proposition \ref{ordersoundness}
(ii) a valuation $u:Var\cup X\mapsto \lbrack 0,1]$ such that $u(\Gamma
_{\varphi ,v}^{1}\cup Th\mathcal{G}_{\square \Diamond })=1$ and $u(\varphi
)<1$. This implies the following relations between $v$ and $u$ that we list
for further use (see Figure 1). Given $\theta _{1},\theta _{2},\theta _{3},$%
\medskip

\noindent \#\textbf{1}. If $v(\Box \theta )>\alpha $ then $u(\theta )=1$
(since then $\theta \in \Gamma _{\varphi ,v}^{1})$

\noindent \#\textbf{2} If $v(\Diamond \theta _{1})\leq v(\square \theta
_{2}) $ then $u(\theta _{1})\leq u(\theta _{2})$ (since then $\theta
_{1}\rightarrow \theta _{2}\in \Gamma _{\varphi ,v}^{1})$

\noindent \#\textbf{3} If $v(\Diamond \theta _{1})<v(\square \theta _{2})$
then $u(\theta _{1})=1$ or $u(\theta _{1})<u(\theta _{2})$ (because then $%
(\theta _{2}\rightarrow \theta _{1})\rightarrow \theta _{1})\in \Gamma
_{\varphi ,v}^{1})$

\noindent \#\textbf{4}. If $u(\theta _{2})<u(\theta _{1})$ then $v(\square
\theta _{2})<v(\Diamond \theta _{1})$ (counter-reciprocal of \textbf{2})

\noindent \#\textbf{5}. If $v(\Box \theta _{2})>0$ then $u(\theta _{2})>0$
(making $\theta _{1}:=\bot $ in \textbf{3} because $u(\bot )=v(\Diamond \bot
)=0$)

\noindent \#\textbf{6}. If $u(\theta _{2})\leq u(\theta _{1})<1,$ then $%
v(\square \theta _{2})\leq v(\Diamond \theta _{1})$ (counter-reciprocal of
\textbf{3}).\newline
\begin{figure}[tbp]
\begin{center}
\epsfysize= 6cm
\epsfig{file=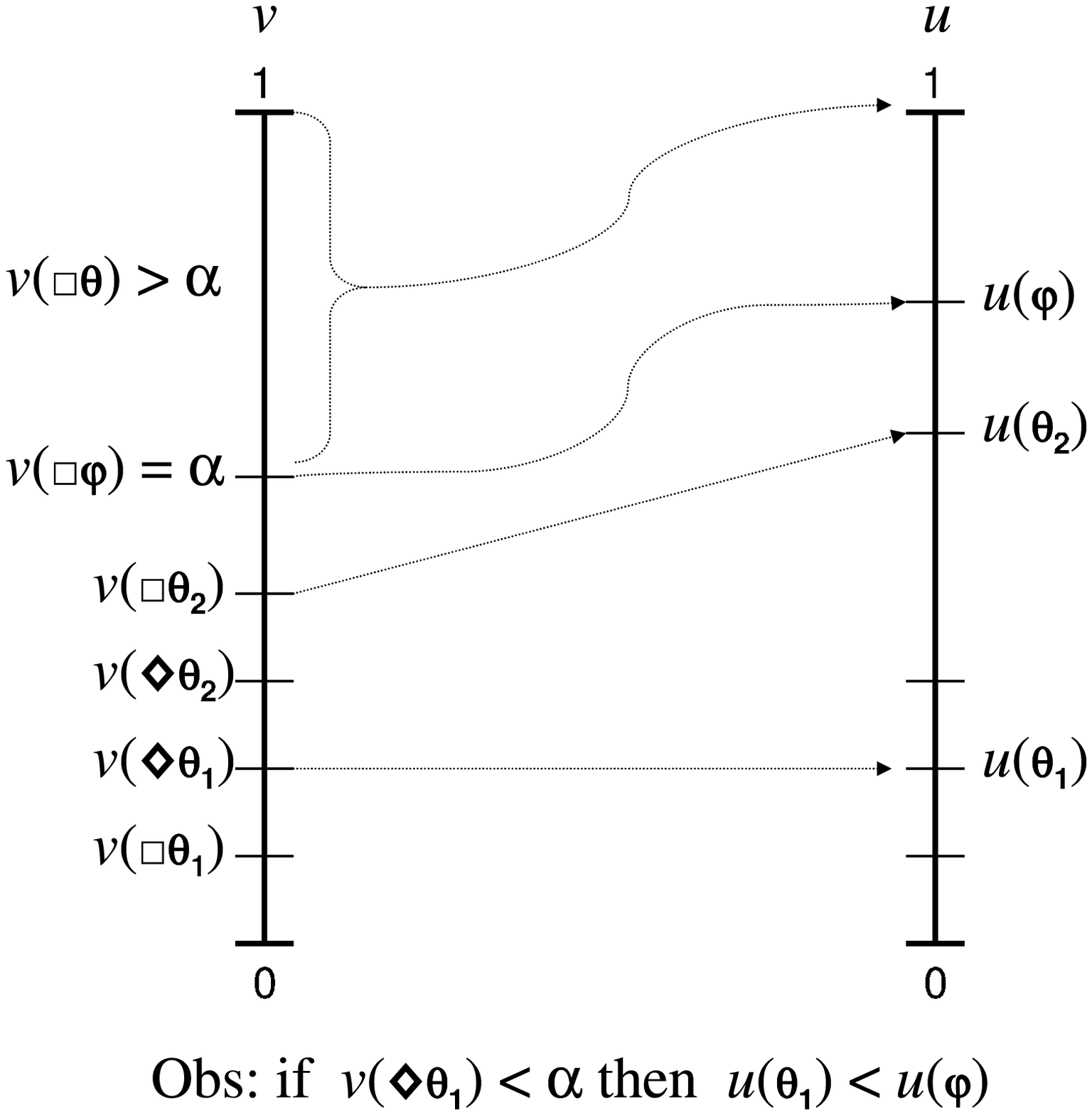, height= 10cm, width=
7cm, angle=270} \vspace{1cm}
\end{center}
\caption{First Translation}
\end{figure}

\medskip

For the next construction we need the finiteness of $F.$ Set $B=\{v(\square
\theta ):\theta \in F\}$ and for each $b\in B$ define%
\begin{equation*}
u_{b}=\min \{u(\theta ):\theta \in F,\text{ }v(\square \theta )=b\},
\end{equation*}%
and then define a strictly descending sequences $b_{0},b_{1},...,b_{N}=0$ in
$B$ as follows:

\medskip

$b_{0}=\alpha $

$b_{i+1}=\max \{b<b_{i}:$ and $u_{b}<u_{b_{i}}\}$\medskip

\noindent Pick formulas $\varphi _{i}\in F$ such that $b_{i}=v(\square
\varphi _{i})$ and $u_{b_{i}}=u(\varphi _{i}).$ By construction the sequence
$u_{b_{0}},u_{b_{1}},...$ is also strictly descending and $%
u_{b_{0}}=u_{\alpha }\leq u(\varphi )<1,$ thus by finiteness of $F$ the
inductive definition ends with some $b_{N}$ (which could be $b_{0}$ in case $%
u_{\alpha }=0).$ To check that $b_{N}=u_{b_{N}}=0,$ assume $b_{N}=v(\square
\varphi _{N})>0,$ then $u_{b_{N}}=u(\varphi _{N})>0=u(\bot )$ by property%
\textbf{\ \#5} above. Since $v(\square \bot )\leq v(\square \varphi _{N})$
then by minimality of $u_{b_{N}}$ we can not have equality, thus $v(\square
\bot )<v(\square \varphi _{N})$ and there exists $b_{N+1}<b_{N},$ a
contradiction. Knowing $v(\square \bot )=b_{N}=0,$ $u_{b_{N}}\leq u(\bot )=0$
by minimality again.

Fix $\varepsilon >0$ such that $\alpha +\varepsilon <1$ and define further
(taking $\min \emptyset =1)$\medskip

$p_{0}=(\alpha +\varepsilon )\cdot \min \{v(\Diamond \theta ):\theta \in F,$
$\alpha <v(\Diamond \theta )\}$

$p_{i}=b_{i+1}\cdot \min \{v(\Diamond \theta ):\theta \in F,$ $%
b_{i}<v(\Diamond \theta )\}$ for $i\geq 1.$\medskip

\noindent We have then $p_{i}>b_{i}$ by finiteness of $F.$ Summing up,
\begin{eqnarray*}
1 &>&\alpha +\varepsilon \geq p_{0}>b_{0}=\alpha \geq p_{1}>b_{1}\geq
....\geq p_{N}>b_{N}=0. \\
1 &>&u_{b_{0}}>u_{b_{1}}>...>u_{b_{N}}=0
\end{eqnarray*}%
Now pick an strictly increasing function $g:[0,1]\mapsto \lbrack 0,1]$ such
that (see Figure 2)\medskip

$g(1)=1$

$g[[u_{\alpha },1)]=[\alpha ,p_{0})$

$g[[u_{b_{i+1}},u_{b_{i}})]=[b_{i+1},p_{i+1})$
\begin{figure}[tbp]
\begin{center}
\epsfysize= 7cm
\epsfig{file=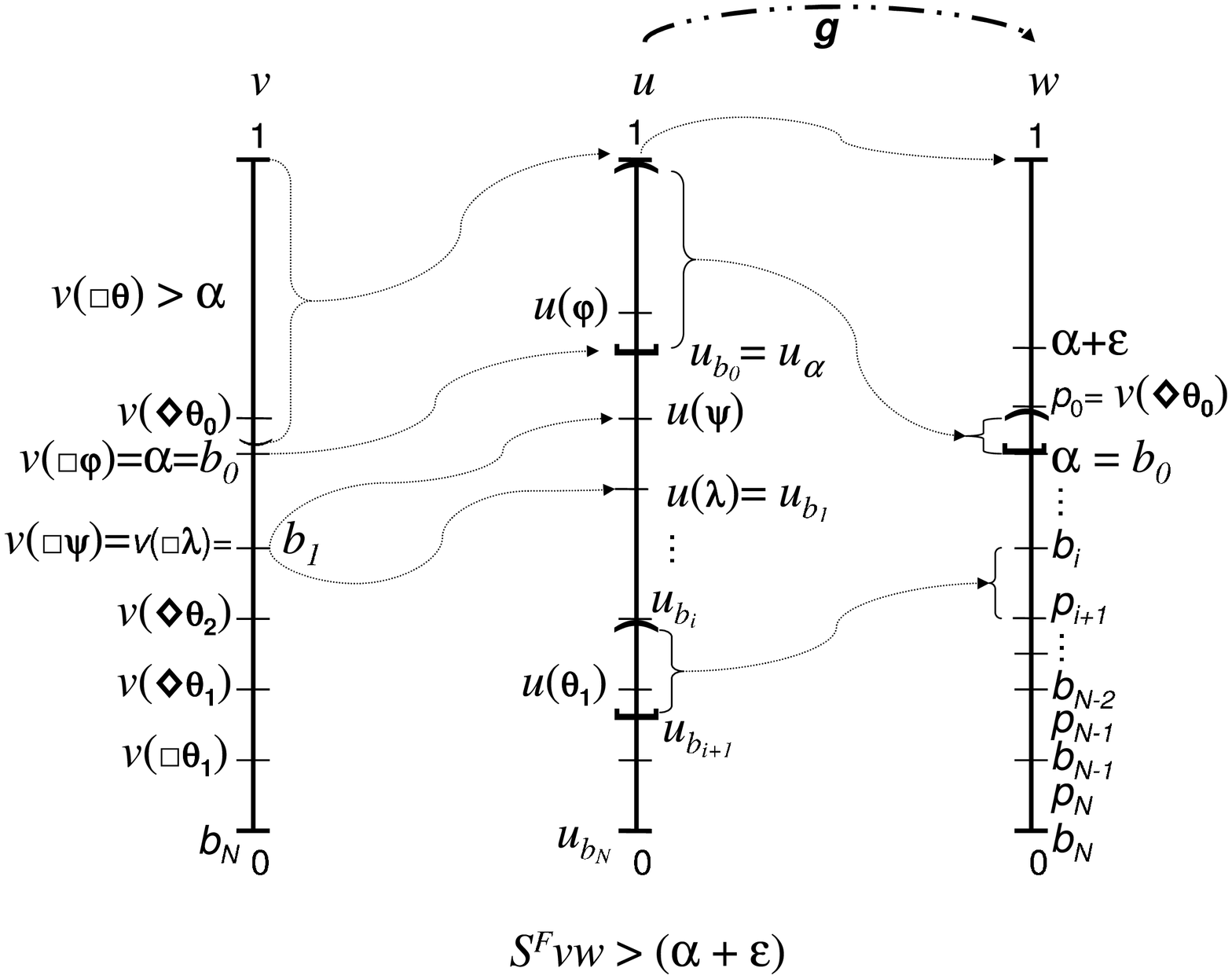, height= 10cm, width=
7cm, angle=270} \vspace{1cm}
\end{center}
\caption{Second Translation}
\end{figure}

\noindent Then the valuation $w=g\circ u$ satisfies $w(Th\mathcal{G}%
_{\square \Diamond })=1,$ and so it belongs to $W.$ Moreover, for any $%
\theta \in F$:

i) If $\ u(\theta )=1$ then $w(\theta )=1$ by definition of $w$;\ hence, $%
v(\square \theta )\leq w(\theta ).$ In addition, $v(\Diamond \theta )\geq
p_{0}$, otherwise $v(\Diamond \theta )\leq \alpha =v(\square \varphi _{0})$
which would imply $u(\theta )\leq $ $u(\varphi _{0})<1$ by \textbf{\#2}, a
contradiction.

ii) If $u(\theta )\in \lbrack u_{b_{i}},u_{b_{i-1}})$ or $u(\theta
)=[u_{b_{0}},1)$ then $v(\square \theta )\leq w(\theta )\leq v(\Diamond
\theta ).$ To see this notice first that $w(\theta )\in \lbrack b_{i},p_{i})$
by definition of $g$. Now, for $i\geq 1$, $b_{i}$ is the maximum $v(\square
\psi )$ with $u(\psi )<u_{b_{i-1}}$ therefore $v(\square \theta )\leq
b_{i}\leq w(\theta ).$ In addition, for $i=0$, $v(\square \theta )\leq
\alpha =b_{0}\leq w(\theta )$ by \textbf{\#1. }Moreover, if $u(\theta
)=u_{b_{i}}=u(\varphi _{i})$ then $w(\theta )=b_{i}=v(\square \varphi
_{i})\leq v(\Diamond \theta )$ by \textbf{\#6}, and if $u(\theta )>u_{b_{i}}$
then $v(\Diamond \theta )>v(\square \varphi _{i})=b_{i}$ by \textbf{\#4};
hence, $v(\Diamond \theta )\geq p_{i}>w(\theta ).$\medskip

From (i,ii), it follows that $\inf_{\theta \in F}\{v(\square \theta
)\Rightarrow w(\theta )\}=1$ and $\inf_{\theta \in F}\{w(\theta )\Rightarrow
v(\Diamond \theta )\}\geq p_{0}.$ Hence $S^{F}vw\geq p_{0},$ and $w(\varphi
)=g(u(\varphi ))<p_{0}\leq \alpha +\varepsilon .$

\medskip

\noindent \textit{Proof of Claim 2}. Assume $v(\Diamond \varphi )=\alpha >0.$%
\begin{equation*}
\begin{array}{ll}
U_{\varphi ,v}^{{}}= & \{\theta :v(\Diamond \theta )<\alpha \} \\
& \cup \{\vartheta _{2}\rightarrow \vartheta _{1}:v(\Diamond \vartheta
_{1})<v(\square \vartheta _{2})\text{ and }v(\Diamond \vartheta _{1})<\alpha
\} \\
& \cup \{(\vartheta _{1}\rightarrow \vartheta _{2})\rightarrow \vartheta
_{1}:v(\Diamond \vartheta _{1})=v(\square \vartheta _{2})\text{ and }%
v(\Diamond \vartheta _{1})<\alpha \}%
\end{array}%
\end{equation*}%
this set is non-empty because $v(\Diamond \bot )=0,$ moreover, for any $\xi
\in U_{\varphi ,v}$ we have $v(\Diamond \xi )<\alpha .$ For the first set by
construction. For the second set of axioms, because $v(\Diamond (\vartheta
_{2}\rightarrow \vartheta _{1}))\leq v(\square \vartheta _{2}\rightarrow
\Diamond \vartheta _{1})=v(\Diamond \vartheta _{1})<\alpha $ by $FS1.$ For
the third, notice that $v(\Diamond ((\vartheta _{1}\rightarrow \vartheta
_{2})\rightarrow \vartheta _{1})))$ $\leq v(\square (\vartheta
_{1}\rightarrow \vartheta _{2})\rightarrow \Diamond \vartheta _{1})$ $\leq
v((\Diamond \vartheta _{1}\rightarrow \square \vartheta _{2})\rightarrow
\Diamond \vartheta _{1})=v(\Diamond \vartheta _{1})<\alpha $ by FS1, FS2.

We claim that for any finite $\{\xi _{1},\ldots ,\xi _{k}\}\subseteq
U_{\varphi ,v}:$%
\begin{equation*}
\varphi \not\vdash _{_{G_{\Box \Diamond }}}\xi _{1}\vee \ldots \vee \xi _{k}
\end{equation*}%
because, on the contrary, $\Diamond \varphi \vdash _{G_{\Box \Diamond
}}\Diamond (\xi _{1}\vee \ldots \vee \xi _{k})\vdash _{G_{\Box \Diamond
}}\Diamond \xi _{1}\vee \ldots \vee \Diamond \xi _{k}$ or
\begin{equation*}
\Diamond \varphi ,Th\mathcal{G}_{\square \Diamond }\vdash \Diamond \xi
_{1}\vee \ldots \vee \Diamond \xi _{k}
\end{equation*}%
and evaluating with $v$ it would give: $\alpha =v(\Diamond \varphi )\leq
\max \{v(\Diamond \xi _{1}),\ldots ,v(\Diamond \xi _{k})\}<\alpha $, absurd.

Therefore, there is a valuation $u$ such that $u(\varphi )=u(T\mathcal{G}%
_{\Box \Diamond })=1$ and $u(\xi )<1$ for each $\xi \in U_{\varphi ,v},$
which has the following consequences for any\ $\theta ,\theta _{1},\theta
_{2}$:\medskip

\noindent \#\#\textbf{1}. If $v(\Diamond \theta )<\alpha $ then $u(\theta
)<1 $ (because then $\theta \in U_{\varphi ,.v})$

\noindent \#\#\textbf{2}. If $v(\Diamond \theta _{1})<v(\Box \theta _{2})$
and $v(\Diamond \theta _{1})<\alpha $ then $u(\theta _{1})<u(\theta _{2})$
(because $\theta _{2}\rightarrow \theta _{1}\in U_{\varphi ,.v})$

\noindent \#\#\textbf{3}. If $v(\Diamond \theta _{1})\leq v(\square \theta
_{2})$ and $v(\Diamond \theta _{1})<\alpha $ then $u(\theta _{1})\leq
u(\theta _{2})$ (because $(\theta _{1}\rightarrow \theta _{2})\rightarrow
\theta _{1}\in U_{\varphi ,v})$

\noindent \#\#\textbf{4 }If\ $u(\theta _{2})=0$ then $v(\Box \theta _{2})=0$
(making $\theta _{1}:=\bot $ in \textbf{2} and taking counter-reciprocal$)$

\noindent \#\#\textbf{5}. If $v(\Diamond \theta _{1})=0$ then $u(\theta )=0$
(making $\theta _{2}:=\bot $ in \textbf{3}, because then $v(\Diamond \theta
_{1})\leq v(\square \bot )$ and $v(\Diamond \theta _{1})<\alpha $).

\medskip

We perform now a construction dual of the one we did in the proof of Claim
1. Let $C=\{v(\Diamond \theta )\leq \alpha :\theta \in F\}$ and define for
each $c\in C$
\begin{equation*}
u_{c}=\max \{u(\theta ):\theta \in F,\text{ }v(\Diamond \theta )=c\}.
\end{equation*}%
Note that $u_{0}=0$ by \textbf{\#\#5} above,\textbf{\ }and $u_{\alpha }=1$
because $u(\varphi )=1.$ Define an ascending sequence $0=c_{0}<c_{1}<....$
in $C$ as follows:\medskip

$c_{0}=v(\Diamond \bot )=0$

$c_{1}=\min \{c\in C:c>c_{0}$ and $u_{c}>u_{c_{0}}\}$

$c_{2}=\min \{c\in C:c>c_{1}$ and $u_{c}>u_{c_{1}}\}$

etc.\medskip

\noindent Choose $\varphi _{i}$ such that $u_{c_{i}}=u(\varphi _{i}),$ $%
c_{i}=v(\Diamond \varphi _{i}),$ clearly, $0=u_{c_{0}}<u_{c_{1}}<....$ By
finiteness of $F$ the sequence of the $c_{i}$ ends necessarily with $%
c_{N}=\alpha ,$ because $c_{i}=v(\Diamond \varphi _{i})<\alpha $ implies $%
u_{c_{i}}=u(\varphi _{i})<1=u_{\alpha }$ by \#\#\textbf{1 }above\textbf{\ }%
and thus the existence of $c_{i+1}\leq \alpha .$ This means also that $%
u_{c_{n}}=1.$

Fix $\varepsilon >0$ such that $\alpha -\varepsilon >c_{N-1},$ and define
further (taking $\max \emptyset =0)$\medskip

$q_{N-1}=\max \{\alpha -\varepsilon ,\max \{v(\square \theta ):v(\square
\theta )<c_{N}\}\}$

$q_{i}=\max \{c_{i},\max \{v(\square \theta ):v(\square \theta
)<c_{i+1}\}\}, $ for $i<N-1$\medskip

\noindent then we have:%
\begin{eqnarray*}
0 &=&c_{0}\leq q_{0}<c_{1}\leq q_{1}<....c_{N-1}\leq \alpha -\varepsilon
\leq q_{N-1}<c_{N}=\alpha \\
0 &=&u_{c_{0}}<u_{c_{1}}<.....<u_{c_{N}}=1
\end{eqnarray*}%
Choose $g:[0,1]\rightarrow \lbrack 0,1]$ to be any strictly increasing
function such that\medskip

$g(0)=0$

$g[(u_{c_{i}},u_{c_{i+1}}]]=(q_{i},c_{i+1}]$ for $i<N-1$

$g[(u_{c_{N-1}},1)]=(q_{N-1},\alpha )$

$g(1)=1$\medskip

\noindent Then $g$ is a Heyting homomorphism and the valuation $w=g\circ v$
satisfies $w(\varphi )=w(T\mathcal{G}_{\Box \Diamond })=1$. Moreover, we
have:\medskip

\noindent i)\ \ If $v(\Diamond \theta )\geq \alpha $ then trivially $%
(w(\theta )\Rightarrow v(\Diamond \theta ))\geq \alpha .$ In particular, $%
(w(\varphi )\Rightarrow v(\Diamond \varphi ))=(1\Rightarrow v(\Diamond
\varphi ))=\alpha .$

\noindent ii)\ If $v(\Diamond \theta )<\alpha $ then $w(\theta )\leq
v(\Diamond \theta ).$ To see this consider cases. First: $u(\theta )\in
(u_{c_{i}},u_{c_{i+1}})$ for some $i$ (recall $u(\theta )<1$ by \#\#\textbf{1%
}) then $w(\theta )\in (q_{i},c_{i+1}]$. As $u(\theta )>u_{c_{i}}$ and $%
c_{i+1}=v(\Diamond \varphi _{i+1})$ is the smallest $v(\Diamond \psi )$ with
$u(\psi )>u_{c_{i}}$ then $v(\Diamond \theta )\geq c_{i+1}\geq w(\theta ).$
Second: $u(\theta )=0$ then $w(\theta )=0$ and $v(\Box \theta )=0$ by
\textbf{\#\#4}

\noindent iii) \ If $v(\square \theta )\geq \alpha $ then $(v(\square \theta
)\Rightarrow w(\theta ))>\alpha -\varepsilon ,$ because $v(\square \theta
)>c_{N-1}=v(\Diamond \varphi _{N-1})$ which implies $u(\theta )>u(\varphi
_{N-1})=u_{c_{N-1}}$ by \textbf{\#\#2,} therefore $w(\theta )>q_{N-1}\geq
\alpha -\varepsilon $ by definition$.$

\noindent iv) $v(\square \theta )<\alpha $ then $v(\square \theta )\leq
w(\theta ).$ To see this notice that $c_{i}\leq v(\square \theta )\leq
q_{i}<c_{i+1}$ for some $i$ and consider cases. First: $v(\square \theta
)=c_{i}=v(\Diamond \varphi _{i})$ then, by\textbf{\ \#\#3, }$%
u_{c_{i}}=u(\varphi _{i})\leq u(\theta ).$ Therefore $c_{i}\leq w(\theta ).$
That is, $v(\square \theta )\leq w(\theta ).$ Second: $c_{i}<v(\square
\theta )$ then, by\textbf{\ \#\#2,\ }$u_{c_{i}}<u(\theta )$ and by
definition $q_{i}\leq w(\theta ),$ which shows again $v(\square \theta )\leq
w(\theta ).$

From (i,ii),\ we have $\inf_{\theta \in F}\{w(\theta )\Rightarrow v(\Diamond
\theta )\}=\alpha $, and from (iii,iv), $\inf_{\theta \in F}\{v(\square
\theta )\Rightarrow w(\theta )\}$ $\geq \alpha -\varepsilon .$ Hence, $%
S^{F}vw\geq \alpha -\varepsilon $ and thus $S^{F}vw\cdot w(\varphi
)=S^{F}vw\cdot \alpha \geq (\alpha -\varepsilon ).$ $\blacksquare $

\begin{lemma}
{\label{JointCompleteness}}\textbf{(Weak completeness) }For any finite
theory $T$ and formula $\varphi $ in $\mathcal{L}_{\square \Diamond }$, $%
T\models _{GK}\varphi $ implies $T\vdash _{\mathcal{G}_{\square \Diamond
}}\varphi .$
\end{lemma}

\proof Assume $T$ is finite and $T\not\vdash _{\mathcal{G}_{\square \Diamond
}}\varphi $ then $T,Th\mathcal{G}_{\square \Diamond }\not\vdash _{\mathcal{G}%
}\varphi $ by Lemma \ref{reduction} and thus there is by Proposition \ref%
{ordersoundness} a Gödel valuation $v:Var\cup X\rightarrow \lbrack 0,1]$
such that $v(\varphi )<v(T)=v(Th\mathcal{G}_{\square \Diamond })=1.$ Let $F$
be the set of subformulas of formulas in $T\cup \{\varphi \}$ (including $%
\bot $ ), then $v\in W.$ the set of worlds of the canonical model $M_{F}$
and by Lemma \ref{equation-joint}, $e^{F}(v,T)=v(T)=1$ and $e^{F}(v,\varphi
)=v(\varphi )<1,$ thus $T\not\models _{GK}\varphi .\blacksquare $\medskip

To prove strong completeness we utilize compactness of first order classical
logic and the following result of Horn:

\begin{lemma}
(\emph{\cite{Horn69}, Lemma 3.7}) Any countable linear order $(P,<)$ may be
embedded in $(\mathbb{Q}\cap \lbrack 0,1],<)$ preserving all joins and meets
existing in $P$.
\end{lemma}

\begin{theorem}
{\label{StrongCompleteness}}\textbf{(Strong completeness) }For any countable
theory $T$ and formula $\varphi $ in $\mathcal{L}_{\square \Diamond }$, $%
T\vdash _{\mathcal{G}_{\square \Diamond }}\varphi $ if and only if $T\models
_{GK}\varphi .$
\end{theorem}

\proof Assume $T$ is countable and $T\nvdash _{\mathcal{G}_{\square \Diamond
}}\varphi $. Consider the first order theory $T^{\ast }$ with two unary
relation symbols\ $W,P,$ binary $<$, constant symbols $0,1$, and $c,$
function symbols $x\circ y,$ $S(x,y),$ and $f_{\theta }(x)$ for each $\theta
\in \mathcal{L}_{\square \Diamond }(V)$ where $V$ is the set of
propositional variables of $T,$ and having for axioms:

\medskip

$\forall x\lnot (Wx\wedge Px)$

$(P,<)$ is a strict linear order with minimum $0$ and maximum $1$

$\forall x\forall y(W(x)\wedge W(y)\rightarrow P(S(x,y)))$

$\forall x\forall y(P(x)\wedge P(y)\rightarrow (x\leq y\wedge x\circ
y=1)\vee (x>y\wedge x\circ y=y))$

$\forall x(W(x)\rightarrow f_{\bot }(x)=0)$

for each $\theta ,\psi \in \mathcal{L}_{\square \Diamond }:$

$\forall x(W(x)\rightarrow P(f_{\theta }(x)))$

$\forall x(W(x)\rightarrow f_{\theta \wedge \psi }(x)=\min \{f_{\theta
}(x),f_{\psi }(x)\})$

$\forall x(W(x)\rightarrow f_{\theta \rightarrow \psi }(x)=(f_{\varphi
}(x)\circ f_{\psi }(x))$

$\forall x(W(x)\rightarrow f_{\square \theta }(x)=\inf_{y}(S(x,y)\circ
f_{\theta }(y))$

$\forall x(W(x)\rightarrow f_{\Diamond \theta }(x)=\sup_{y}(\min
\{S(x,y),f_{\theta }(y)\})$

$W(c)\wedge (f_{\varphi }(c)<1)$

for each $\theta \in T:f_{\theta }(c)=1$

\medskip

\noindent For each finite part $t$ of this theory let $\Sigma _{t}=\{\theta
:f_{\theta }\in t\}.$ Since $\Sigma _{t}\cap T\nvdash _{\mathcal{G}_{\square
\Diamond }}\varphi \ $by hypothesis then by weak completeness there is a
GK-model $\,M_{\Sigma }=(W_{\Sigma },S_{\Sigma },e_{\Sigma })$ and $a\in
W_{\Sigma }$ such that $e_{\Sigma }(a,\theta )=1$ for each $\theta \in
\Sigma _{t}\cap T$ and $e_{\Sigma }(a,\varphi )<1.$ Therefore the first
order structure $(W_{\Sigma }\sqcup \lbrack 0,1],W_{\Sigma
},[0,1],<,0,1,a,\Rightarrow ,S_{\Sigma },f_{\theta })_{\theta \in \mathcal{L}%
_{\square \Diamond }},$ with $f_{\theta }:W_{\Sigma }\rightarrow \lbrack
0,1] $ defined as $f_{\theta }(x)=e_{\Sigma }(x,\theta ),$ is clearly a
model of $t.$ By compactness of first order logic and the downward Löwenheim
theorem $T^{\ast }$ has a countable model $M^{\ast }=(B,W,P,<,0,1,a,\circ
,S,f_{\theta })_{\theta \in \mathcal{L}_{\square \Diamond }}.$ Using Horn's
lemma \cite{Horn69}, $(P,<)$ may be embedded in $(\mathbb{Q}\cap \lbrack
0,1],<)$ preserving $0,1,$ and all suprema and infima existing in $P$;
therefore, we may assume without loss of generality that the ranges of the
functions $S$ and $f_{\theta }$ are contained in $[0,1].$ Then, it is
straightforward to verify that $M=(W,S,e),$ where $e(x,\theta )=f_{\theta
}(x)$ for all $x\in W,$ is a GK-model such that $M\models _{a}T,$ and $%
M\not\models _{a}\varphi ,$ that is, $T\not\models _{GK}\varphi .$ The rest
follows from Theorem \ref{strong soundness}. $\blacksquare $

\medskip

We can not expect a similar result for uncountable theories by the
observation after Proposition \ref{ordersoundness}. If we allow non-standard
values, for example $HK$-models where $H$ is an adequate ultrapower of
[0,1], we may obtain strong completeness for entailment of uncountable
theories up to certain cardinality. However,

\begin{proposition}
{\label{uncountable}}There is no single linearly ordered Heyting algebra $H$
giving strong completeness with respect to HK models for theories of
arbitrary power, even in Gödel-Dummet logic.
\end{proposition}

\proof Assume otherwise, then $H$ would be infinite (by the old Gödel
argument). Let $\kappa $ be a cardinal greater than $|H|$ and consider the
theory $T=\{(p_{\beta }\rightarrow p_{\alpha })\rightarrow q:\alpha <\beta
<\kappa \},$ then $T\models _{HK}q$ because $v(T)=1$ with $v(q)<1$ would
imply $v(p_{\beta }\rightarrow p_{\alpha })<1$ and thus $v(p_{\alpha
})<v(p_{\beta })$ for $\alpha <\beta <\kappa ,$ which is impossible by
cardinality considerations. On the other hand, $T\nvdash _{\mathcal{G}_{\Box
\Diamond }}q,$ otherwise we would have $\Delta \vdash _{\mathcal{G}_{\Box
\Diamond }}q$ and thus $\Delta \models _{HK}q,$ for some finite set $\Delta
=\{(p_{\alpha _{i+1}}\rightarrow p_{\alpha _{i}})\rightarrow q:1\leq i<n\},$
which is impossible because any valuation $v(q)=h,$ $v(p_{\alpha
_{i}})=h_{i} $ where $h_{1}<h_{2}<...<h_{n+1}<h<1$ makes $v((p_{\alpha
_{i+1}}\rightarrow p_{\alpha _{i}})\rightarrow q)=1$ for $1\leq i<n.$ $%
\square $

\section{Optimal models, modal axioms}

To extend the completeness theorem to the [0,1]-valued analogues of the
classical bi-modal systems $T,$ $S4$, $S5$ we introduce a particular kind of
GK-model, their advantage being that the many-valued counterpart of
classical structural properties of frames may be characterized in them by
the validity of the corresponding classical schemes.

Given a GK-model $M=(W,S,e),$ define a new accessibility relation $S^{+}xy=$
$S_{\square }xy\cdot $ $S_{\Diamond }xy,$ where $S_{\square
}xy=\inf_{\varphi \in \mathcal{L}_{\square \Diamond }}\{e(x,\Box \varphi
)\Rightarrow e(y,\varphi )\},$\ and\ $S_{\Diamond }xy=\inf_{\varphi \in
\mathcal{L}_{\square \Diamond }}\{e(y,\varphi )\Rightarrow e(x,\Diamond
\varphi )\},$ and call $M$ \emph{optimal} if $S^{+}=S.$

The following lemma shows that any model is equivalent to an optimal one.

\begin{lemma}
{\label{optimal}}$(W,S^{+},e)$ is optimal and if $e^{+}$ is the extension of
$e$ in this model then $e^{+}(x,\varphi )=e(x,\varphi )$ for any $\varphi
\in \mathcal{L}_{\square \Diamond }.$
\end{lemma}

\proof The first claim follows from the second which is proven by a
straightforward induction on formulas. The only non trivial step is that of
the modal connectives. Notice first that $Sxy\leq S^{+}xy,$ because $%
e(x,\Box \varphi )\leq (Sxy\Rightarrow e(y,\varphi ))$ and $Sxy\cdot
e(y,\varphi )\leq e(x,\Diamond \varphi )$ for any $\varphi $, thus $Sxy\leq
(e(x,\Box \varphi )\Rightarrow e(y,\varphi )),(e(y,\varphi )\Rightarrow
e(x,\Diamond \varphi )).$ Now, assume $e^{+}(y,\varphi )=e(y,\varphi )$ for
all $y$ then by the first observation and the induction hypothesis, $%
e^{+}(x,\square \varphi )=\inf_{y}\{S^{+}xy\Rightarrow e^{+}(y,\varphi
)\}\leq \inf_{y}\{Sxy\Rightarrow e(y,\varphi )\}=e(x,\square \varphi )$. But
$S^{+}xy\leq (e(x,\square \varphi )\Rightarrow e(y,\varphi ))$ by definition
of $S^{+}$ and thus $e(x,\square \varphi )\leq (S^{+}xy\Rightarrow
e(y,\varphi ))=(S^{+}xy\Rightarrow e^{+}(y,\varphi ))$ which yields $%
e(x,\square \varphi )\leq e^{+}(x,\square \varphi ).$ Similarly, by the
induction hypothesis and the first observation, $e^{+}(x,\Diamond \varphi
)=\sup_{y}\{S^{+}xy\cdot e(y,\varphi )\}\geq \sup_{y}\{Sxy\cdot e(y,\varphi
)\}=e(x,\Diamond \varphi ),$ and by definition $S^{+}xy\leq (e(y,\varphi
)\Rightarrow e(x,\Diamond \varphi ))$ and thus $S^{+}xy\cdot e^{+}(y,\varphi
)=S^{+}xy\cdot e(y,\varphi )\leq e(x,\Diamond \varphi )$ which yields $%
e^{+}(x,\Diamond \varphi )\leq e(x,\Diamond \varphi )$. $\blacksquare $%
\medskip

Call a GK-frame $\mathcal{M}$ = $\langle W,S\rangle $ \emph{reflexive} if $%
Sxx=1$ for all $x\in W$, \emph{transitive} if $Sxy\cdot Syz\leq Sxz$ for all
$x,y,z,$ and \emph{symmetric} if $Sxy=Syx$ for all $x,y\in W.$ Let $Ref$, $%
Trans,$ and $Symm$ denote, respectively, the classes of GK- models over
frames satisfying each one of the above properties. These are the fuzzy
versions of the corresponding classical properties of frames characterized
by the following pairs of modal axioms:%
\begin{equation}
\begin{array}{lllllllll}
\text{T}_{\square }. & \Box \varphi \rightarrow \varphi &  &  & \text{T}%
_{\Diamond }. & \varphi \rightarrow \Diamond \varphi &  &  & \text{%
reflexivity} \\
\text{4}_{\square }\text{.} & \Box \varphi \rightarrow \Box \Box \varphi &
&  & \text{4}_{\Diamond }. & \Diamond \Diamond \varphi \rightarrow \Diamond
\varphi &  &  & \text{transitivity} \\
\text{M}_{1}. & \varphi \rightarrow \square \Diamond \varphi &  &  & \text{M}%
_{2}. & \Diamond \square \varphi \rightarrow \varphi &  &  & \text{symmetry}%
\end{array}
\label{table}
\end{equation}%
We will see that these axioms characterize also the fuzzy versions in
optimal models.

\begin{lemma}
\label{similarity}i)$\ $T$_{\square }$ and T$_{\Diamond }$ are valid in $Ref$%
. ii)$\ 4_{\square }$ and $4_{\Diamond }$ are valid in $Trans.$ iii)$\ $M$%
_{1}$ and M$_{2}$ are\textbf{\ }valid in $Symm.$
\end{lemma}

\proof i)\ In reflexive models, $e(x,\Box \varphi )\leq (Sxx\Rightarrow
e(x,\varphi ))=e(x,\varphi )$ and $e(x,\Diamond \varphi )$ $\geq Sxx\cdot
e(x,\varphi )$ $=e(x,\varphi )$ for any $x.$ Thus $e(x,\Box \varphi
\rightarrow \varphi )=1=e(x,\varphi \rightarrow \Diamond \varphi ).$

ii)\ In transitive models $e(x,\Box \varphi )\cdot Sxy\cdot Syz\leq \lbrack
(Sxz\Rightarrow e(z,\varphi ))\cdot Sxz]\leq e(z,\varphi )$ for all $x,y,z.$
Hence, $e(x,\Box \varphi )\cdot Sxy\leq (Syz\Rightarrow e(z,\varphi ))$ and
thus $e(x,\Box \varphi )\cdot Sxy\leq e(y,\Box \varphi );$ therefore, $%
e(x,\Box \varphi )\leq (Sxy\Rightarrow e(y,\Box \varphi ))$ for all $y$ and
thus $e(x,\Box \varphi )\leq e(x,\Box \Box \varphi )$ which yields $%
4_{\square }$. Also $Sxy\cdot Syz\cdot e(z,\varphi )$ $\leq Sxz\cdot
e(z,\varphi )\leq e(x,\Diamond \varphi ).$ Hence, $Syz\cdot e(z,\varphi
)\leq (Sxy\Rightarrow e(x,\Diamond \varphi ))$ and thus $e(x,\Diamond
\varphi )\leq (Sxy\Rightarrow e(x,\Diamond \varphi ));$ therefore, $Sxy\cdot
e(x,\Diamond \varphi )\leq e(x,\Diamond \varphi ))$ for all $y$ and thus $%
e(x,\Diamond \Diamond \varphi )\leq e(x,\Diamond \varphi )$ which gives $%
4_{\Diamond }.$ iii) In symmetric models, $Sxy\cdot e(x,\varphi )=Syx\cdot
e(x,\varphi )\leq e(y,\Diamond \varphi )$ for all $x,y$, then $e(x,\varphi
)\leq (Sxy\Rightarrow e(y,\Diamond \varphi ))$ and thus $e(x,\varphi )\leq
e(y,\square \Diamond \varphi ))$ which is M$_{1}.$ Moreover, $e(y,\square
\varphi )\leq (Syx\Rightarrow e(x,\varphi )),$ thus $Sxy\cdot e(y,\square
\varphi )=Syx\cdot e(y,\square \varphi )\leq e(x,\varphi )$ and thus $%
e(x,\Diamond \square \varphi )\leq e(x,\varphi )$ which is M$_{2}.$ $%
\blacksquare $

\begin{proposition}
\label{character}Let $M$ be an optimal GK-model, then \ i) It is reflexive
if and only if it validates the schemes T$_{\square }+$T$_{\Diamond }.$ ii)$%
\ $It is transitive if and only if it validates $4_{\square }+4_{\Diamond }.$
iii) It is symmetric if and only if it validates M$_{1}+$M$_{2}$.
\end{proposition}

\proof i)\ By optimality, $Sxx=\inf_{\varphi }\{e(x,\square \varphi
\rightarrow \varphi )\}\cdot \inf_{\varphi }\{e(x,\varphi \rightarrow
\Diamond \varphi )\}=1$ if T$_{\square },$T$_{\Diamond }$ hold in $M$. ii)\
By definition,

$%
\begin{array}{ll}
S_{\square }xy\cdot S_{\square }yz & \leq (e(x,\square \square \varphi
)\Rightarrow e(y,\square \varphi ))\cdot (e(y,\square \varphi )\Rightarrow
e(z,\varphi )) \\
& \leq (e(x,\square \square \varphi )\Rightarrow e(z,\varphi ))\leq
(e(x,\square \varphi )\Rightarrow e(z,\varphi )),%
\end{array}%
$

\noindent the last inequality holding by $4_{\square }$. Similarly

$%
\begin{array}{ll}
S_{\Diamond }xy\cdot S_{\Diamond }yz & \leq (e(y,\Diamond \varphi
)\Rightarrow e(x,\Diamond \Diamond \varphi ))\cdot (e(z,\varphi )\Rightarrow
e(y,\Diamond \varphi )) \\
& \leq (e(z,\varphi )\Rightarrow e(x,\Diamond \Diamond \varphi ))\leq
(e(z,\varphi )\Rightarrow e(x,\Diamond \varphi )%
\end{array}%
$

\noindent the last inequality holding by $4_{\Diamond },$ and taking meet
over $\varphi $ in the right we get transitivity. iii)\ since $S_{\square
}xy\leq (e(x,\square \Diamond \varphi )\Rightarrow e(y,\Diamond \varphi ))$ $%
\leq (e(x,\varphi )\Rightarrow e(y,\Diamond \varphi ))$ by M$_{1},$ then
taking meet over $\varphi ,$ we obtain $S_{\square }xy\leq S_{\Diamond }yx.$
Similarly, $S_{\Diamond }yx$ $\leq (e(x,\square \varphi )\Rightarrow
e(y,\Diamond \square \varphi ))$ $\leq (e(x,\square \varphi )\Rightarrow
e(y,\varphi ))$ by M$_{2},$ and then $S_{\Diamond }yx\leq S_{\square }xy.$
From this, $S_{\Diamond }xy=S_{\square }yx,$ and\ thus $Sxy=Syx.$ $%
\blacksquare $

\bigskip

\textbf{Remark}. The notion of optimality and all the results in this
section make sense and hold for $HK$-models, for any complete Heyting
algebra $H$.

\section{Gödel analogues of classical bi-modal systems}

Lemma \ref{optimal} in conjunction with Proposition \ref{character} implies
strong completeness of any combination of the axiom pairs in table \ref%
{table} with respect to GK-frames satisfying the associated structural
properties. In particular, we have for the analogues of the classical modal
systems T, S4 and S5:\medskip

$%
\begin{array}{lll}
\mathcal{G}T_{\square \Diamond } & :=\mathcal{G}_{\square \Diamond
}+T_{\square }+T_{\Diamond } &  \\
\mathcal{G}S4_{\square \Diamond } & :=\mathcal{G}T_{\square \Diamond
}+4_{\square }+4_{\Diamond } &  \\
\mathcal{G}S5_{\square \Diamond } & :=\mathcal{G}S4_{\square \Diamond
}+M1+M2. &
\end{array}%
$

\begin{theorem}
\label{sistemascompletos}i)$\ \mathcal{G}T_{\square \Diamond }$ is strongly
complete for $\models _{Ref}.$ ii)$\ \mathcal{G}S4_{\square \Diamond }$ is
strongly complete for $\models _{Ref\cap Trans}.$ iii)$\ \mathcal{G}%
S5_{\square \Diamond }$ is strongly complete for $\models _{Ref\cap
Trans\cap Symm}$.
\end{theorem}

\proof i)\ If $T\models _{Ref}\varphi $ then $T\models _{\text{Optimal}\cap
Ref}\varphi $, thus $T+\{T_{\square },T_{\Diamond }\}\models _{\text{Optimal
}}\varphi $ by Proposition \ref{character}, and $T+\{T_{\square
},T_{\Diamond }\}\models _{\text{GK}}\varphi $ by Lemma \ref{optimal}%
.Therefore, $T+\{T_{\square },T_{\Diamond }\}\vdash _{\mathcal{G}_{\square
\Diamond }}\varphi $, which implies $T\vdash _{\mathcal{G}T_{\square
\Diamond }}\varphi .$ Here, we have used $\{...\}$ to denote the set of all
instances of the schemes within the brackets. The proofs of (ii)\ and (iii)
are similar. $\square \medskip $

After some calculation $\mathcal{G}S5_{\square \Diamond }$ may be seen
deductively equivalent to Prior's and Bull system $MIPC$ \cite{Prior}, \cite%
{Bull} plus the prelinearity axiom$\ (\varphi \rightarrow \psi )\vee (\psi
\rightarrow \varphi ).$ This system presents some interesting features with
respect to its Gödel-Kripke semantics, given by GK-frames $(W,S)$ where $S$
is a fuzzy equivalence relation. Although the uni-modal fragments of $%
\mathcal{G}T_{\square \Diamond }$ and $\mathcal{G}S4_{\square \Diamond }$
are axiomatizable by the double negation shift axioms and the proper axioms
in \ref{table}, as shown \cite{CaicedoRodriguez}, the intrinsic
axiomatization of the uni-modal fragments of $\mathcal{G}S5_{\square
\Diamond }$ remains open. Moreover, in distinction of the other modal
systems here considered, the $\square $-fragment of $\mathcal{G}S5_{\square
\Diamond }$ is not characterized by models with crisp accessibility
relation, as the following example illustrates.$\medskip $

\noindent \textbf{Example} The formula $\square (\square \varphi \vee \psi
)\rightarrow (\square \varphi \vee \square \psi )$ is not a theorem of $%
\mathcal{G}S5_{\square }$ but it is valid in any accessibility-crisp model
of $\mathcal{G}S5_{\square }$. The first claim is granted by the following
counter model:\newline

\begin{center}
\epsfig{file=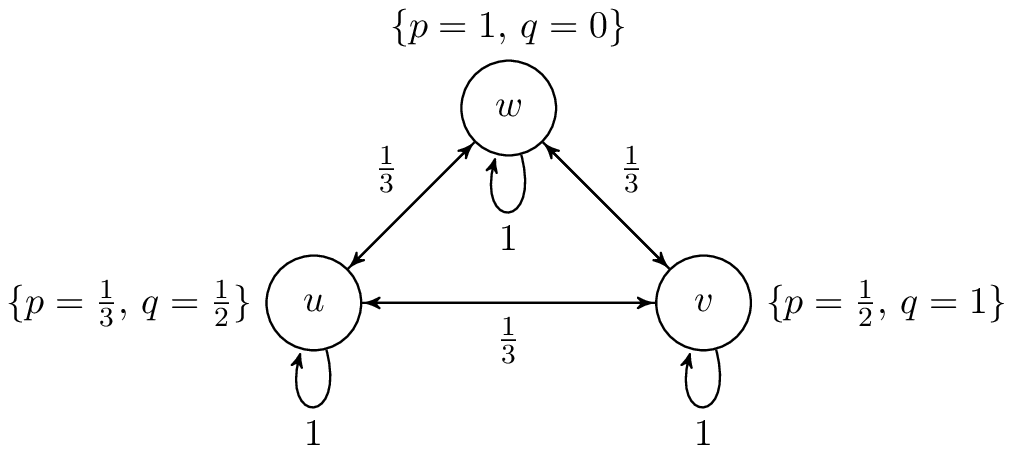, height= 5cm, width= 10cm, angle=0}
\end{center}

\noindent in which the reader may verify that $e(u,\square (\square p\vee q)=%
\frac{1}{2}$ and $e(u,\square p\vee \square q)=\frac{1}{3}$. To prove the
second claim notice that if $(W,S,e)\in Ref\cap Trans\cap Symm$ has crisp $%
S, $ this defines a classical equivalence relation $\sim $\ and thus $%
e(x,\square \varphi )=\inf_{y}\{Sxy\Rightarrow e(y,\varphi )\}=\inf_{y\sim
x}\{e(y,\varphi )\}$ for any formula $\varphi $. Therefore, $e(x,\square
(\square \varphi \vee \psi ))=$ $\inf_{y\sim x}\{\inf_{z\sim y}e(z,\varphi
)\curlyvee e(y,\varphi )\},$ but $\{z:z\sim y\}=\{z:z\sim x\}$ and so $%
\alpha _{y}=$ $\inf_{z\sim y}e(z,\varphi )$ is independent of $y$ for any $%
y\sim x$; hence, $e(x,\square (\square \varphi \vee \psi ))=$ $\inf_{y\sim
x}\{\alpha \curlyvee e(y,\varphi )\}=\alpha \curlyvee \inf_{y\sim
x}\{e(y,\varphi )\}=$ $e(x,\square \varphi \vee \square \psi )$ by
properties of $[0,1]$.

\section{The algebraic connection}

As an algebrizable deductive logic, $\mathcal{G}_{\square \Diamond }$ has a
unique algebraic semantics given by the variety of \emph{bi-modal Gödel
algebras} $A=(G,I,K)$ where $G$ is a Gödel algebra and $I$ and $K$ are unary
operations in $G$ satisfying the identities:\medskip

$%
\begin{array}{lllll}
I(a\cdot b)=Ia\cdot Ib &  &  & K(a\curlyvee b)=Ka\curlyvee Kb &  \\
I1=1 &  &  & K0=0 &  \\
Ka\rightarrow Ib\leq I(a\rightarrow b) &  &  & K(a\rightarrow b)\leq
Ia\rightarrow Kb &
\end{array}%
$\medskip

\noindent This means that $\mathcal{G}_{\square \Diamond }$ is complete with
respect to valuations $v:Var\rightarrow A$ in these algebras, when they are
extend to $\mathcal{L}_{\square \Diamond }$ interpreting $\square $ and $%
\Diamond $ by $I$ and $K$, respectively.

Similarly, $\mathcal{G}T_{\square \Diamond }$, $\mathcal{G}S4_{\square
\Diamond },$ and $\mathcal{G}S5_{\square \Diamond }$ have for algebraic
semantic the subvarieties of bi-modal Gödel algebras determined by the
corresponding pairs of identities in the following table:

\begin{equation}
\begin{array}{llllll}
Ia\leq a &  & a\leq Ka &  &  & \text{reflexivity} \\
Ia=IIa &  & Ka=KKa &  &  & \text{transitivity} \\
a\leq IKa &  & KIa\leq a &  &  & \text{symmetry}%
\end{array}
\label{Table2}
\end{equation}

\noindent Notice that the algebraic models of $\mathcal{G}S4_{\square
\Diamond }$ are just the bi-topological pseudo-Boolean algebras of Ono \cite%
{Ono77} with linear underlying Heyting algebra, and the algebraic models of $%
\mathcal{G}S5_{\square \Diamond }$ are the the monadic Heyting algebras of
Monteiro and Varsavsky \cite{MonVars57}, utilized later by Bull and Fischer
Servi to interpret MIPC, with a Gödel basis. It is proper to call them \emph{%
monadic Gödel algebras}.\medskip

\noindent \textbf{Example.} As we have noticed, there is no finite
counter-model for the formula $\Box \lnot \lnot p\rightarrow \lnot \lnot
\Box p$ in Gödel-Kripke semantics. However, the algebra $A=(\{0,a,1\},I,K)$
where $\{0<a<1\}$ is the three elements Gödel algebra and \ $I1=1,$ $%
Ia=I0=0, $ $K1=Ka=1,$ $K0=0$ is a bi-modal Gödel algebra (actually a monadic
Heyting algebra)\ providing a finite counterexample to the validity of the
formula by means of the valuation $v(p)=a$, as the reader may verify.\medskip

We may associate to each Gödel-Kripke frame $\mathcal{F}=(W,S)$ a bi-modal Gö%
del algebra $[0,1]^{\mathcal{F}}=([0,1]^{W},I^{\mathcal{F}},K^{\mathcal{F}})$
where $[0,1]^{W}$ is the product Gödel algebra, and for each map $f\in
\lbrack 0,1]^{W}:$%
\begin{eqnarray*}
I^{\mathcal{F}}(f)(w) &=&\inf_{w^{\prime }\in W}(Sww^{\prime }\Rightarrow
f(w^{\prime })) \\
K^{\mathcal{F}}(f)(w) &=&\sup_{w^{\prime }\in W}(Sww^{\prime }\cdot
f(w^{\prime }))
\end{eqnarray*}

\begin{theorem}
\label{complexalgebras}$[0,1]^{\mathcal{F}}\ $is a bi-modal Gödel algebra,
and there is a one to one correspondence between Gödel Kripke models over $%
\mathcal{F}$, and valuations $v:Var\rightarrow \lbrack 0,1]^{\mathcal{F}}$
given by the adjunction:%
\begin{equation*}
Var\times W\overset{e}{\rightarrow }[0,1]\text{ \ }\leftrightarrow \text{ \ }%
Var\overset{v_{e}}{\rightarrow }[0,1]^{W},\text{ \ }v_{e}(p)=e(-,p)
\end{equation*}%
so that for any formula $\varphi $, $v_{e}(\varphi )=e(_{-},\varphi ).$

Moreover, the transformation $\mathcal{F}\longmapsto \lbrack 0,1]^{\mathcal{F%
}}$ preserves reflexivity, transitivity and symmetry. Thus, it send Gö%
del-Kripke frames for $\mathcal{G}T_{\square \Diamond }$, $\mathcal{G}%
S4_{\square \Diamond },$ and $\mathcal{G}S5_{\square \Diamond }$ into
algebraic models for the same logics.
\end{theorem}

\proof The verification of the identities that $I^{\mathcal{F}},$ $K^{%
\mathcal{F}}$ must satisfy in each case is routine and the induction in
formulas showing $v_{e}(\varphi )(w)=e(w,\varphi )$ is straightforward. $%
\square \medskip $

Call an algebra of the form $[0,1]^{\mathcal{F}}$ a \emph{Gödel complex
algebra}. Going from algebras to GK-models seems more difficult. However,
utilizing a refinement of our strong completeness theorem for Gödel-Kripke
semantics we may associate to each countable bi-modal Gödel algebra $A$ a
GK-frame $\mathcal{F}_{A}$ such that $A\ $may be embedded in the associated
algebra $[0,1]^{\mathcal{F}_{A}},$ and to each algebraic valuation $\eta $
in $A$ a GK-model over $\mathcal{F}_{A}$ validating the same formulas as $%
\eta $.

Call a theory $T\subseteq \mathcal{L}_{\square \Diamond }$ \emph{normal }if $%
T$ $\vdash _{\mathcal{G}_{\square \Diamond }}\theta $ implies $T\vdash _{%
\mathcal{G}_{\square \Diamond }}\square \theta $ and $T$ $\vdash _{\mathcal{G%
}_{\square \Diamond }}\theta \rightarrow \rho $ implies $T$ $\vdash _{%
\mathcal{G}_{\square \Diamond }}\Diamond \theta \rightarrow \Diamond \rho $.

It $T$ is normal, then for each finite fragment $F$ the proof of the Lemma %
\ref{equation-joint} goes through for the submodel $%
M_{F}^{T}=(W^{T},S^{F},e^{F})$ of the canonical model where $W^{T}=\{v\in
W:v(T)=1\}$. Hence, if $\Sigma $ is a finite subset of $T$ such that $\Sigma
\not\vdash _{\mathcal{G}_{\square \Diamond }}\varphi $ there is a canonical
model $M_{F}^{T}$ such that $e^{F}(v,\Sigma )=1$ and $e^{F}(v,\varphi )<1$
(take $F\supseteq \Sigma \cup \{\varphi \}).$

\begin{lemma}
If $T$ is a countable normal theory there is GK-model $M_{T}$ such that $%
T\vdash _{\mathcal{G}_{\square \Diamond }}\mathcal{\varphi }$ if and only if
$M_{T}\models \mathcal{\varphi }$.
\end{lemma}

\proof From the previous observation, and utilizing a compactness argument
as in the proof of Theorem \ref{StrongCompleteness} we may pick for each $%
\varphi $ such that $T\not\vdash _{\mathcal{G}_{\square \Diamond }}\varphi $
a model $M_{\varphi }=\langle W_{\varphi },S_{\varphi },e_{\varphi }\rangle $
such that $e_{\varphi }(w,T)=1$ for all $w$ and $e(w_{\varphi },\varphi )<1$%
. Define $M_{T}=(W,S,e)$ where $W=\amalg _{\varphi }W_{\varphi },$ $%
Sww^{\prime }=S_{\varphi }ww^{\prime }$ if\ $w,w^{\prime }\in W_{\varphi }$
and 0 otherwise, and $e(w,p)=e_{\varphi }(w,p)$ for $w\in W_{\varphi }.$ It
is easily verified by induction in the complexity of $\theta $ that $%
e(w,\theta )=e_{\varphi }(w,\theta )$ for any $w\in W_{\varphi }.$ Since $%
M_{T}\models T$ then $T\vdash _{\mathcal{G}_{\square \Diamond }}\mathcal{%
\varphi }$ implies $M_{T}\models \mathcal{\varphi }$; reciprocally, if $%
T\not\vdash _{\mathcal{G}_{\square \Diamond }}\mathcal{\varphi }$ then $%
e(w_{\varphi },\varphi )=e_{\varphi }(w_{\varphi },\varphi )<1$ by
construction and thus $M_{T}\not\models \mathcal{\varphi }$. $\square $

\begin{theorem}
For any countable bi-modal Gödel algebra $A$ there is Gödel frame $\mathcal{F%
}_{A}=(W,S)$ such that:

i)\ $A$ is embeddable in the Gödel complex algebra $[0,1]^{\mathcal{F}_{A}}.$

ii)\ For any valuation $v:Var\rightarrow A$ there is a $e_{v}:W\times
Var\rightarrow \lbrack 0,1]$ such that $v(\varphi )=1$ if and only if $%
(W,S,e_{v})\models \varphi .$
\end{theorem}

\proof Fix a valuation $\eta $ into $A$ with onto extension $\eta :\mathcal{L%
}_{\square \Diamond }\rightarrow A$ and let $T=\{\theta :\eta (\theta )=1\},$
then $T$ is normal and for the model $M_{T}=(W,S,e)$ of the previous theorem
we have $\eta (\varphi )=1$ if and only if $e(w,\varphi )=1$ for all $w\in
W. $

(\textit{i}) By Theorem \ref{complexalgebras},\ $e$ induces a bi-modal Gödel
valuation $v_{e}:Var\rightarrow \lbrack 0,1]^{(W,S)},$ $v_{e}(p)=e(-,p)$
such that $v_{e}(\varphi )=e(-,\varphi )=\mathbf{1}\in \lbrack 0,1]^{W}$ if
an only if $\eta (\varphi )=1$ by the observation above. This means that the
extension $v_{e}:\mathcal{L}_{\square \Diamond }\rightarrow \lbrack
0,1]^{(W,S)}$ factors injectively through $\eta $; that is, $v_{e}=\delta
\circ \eta $ for an injective homomorphism of bi-modal Gödel algebras $%
\delta :A\rightarrow \lbrack 0,1]^{(W,S)}$,which shows (\textit{i}). To see
(ii)\ let $v:Var\rightarrow A$, then $\delta \circ v$ is a valuation into $%
[0,1]^{(W,S)}$ which induces, by Theorem \ref{complexalgebras}, a
GK-valuation $e_{v}:W\times Var\rightarrow \lbrack 0,1]$ such that $%
e_{v}(w,\varphi )=\delta (v(\varphi ))(w).$ As $\delta $ is one to one we
have that $v(\varphi )=1$ if and only if $\delta (v(\varphi ))=\mathbf{1}\in
\lbrack 0,1]^{W}$; that is, $e_{v}(w,\varphi )=1$ for all $w,$ which means $%
(W,S,e_{\mu v})\models \varphi .$\ $\square $

$\medskip $

Applying part (\textit{i}) of the previous theorem to the free bi-modal
algebra of countable rank we obtain:

\begin{theorem}
The complex algebras generate the variety of bi-modal Gödel algebras. A
similar result holds for the subvarieties determined by any combination of
identities in (\ref{Table2}).$\medskip $
\end{theorem}

\end{document}